\documentclass[reqno,11pt,a4paper]{amsart}
\usepackage{latexsym,amssymb,amsfonts}
\usepackage{graphicx,color,graphicx}
\usepackage{mathrsfs}
\usepackage{lscape}
\newtheorem{lemma}{\bf Lemma}[section]

\newtheorem{remark}{Remark}

\newtheorem{theorem}{\bf Theorem}[section]

\theoremstyle{plain} \numberwithin{equation}{section}
\theoremstyle{plain} \theoremstyle{definition}

\theoremstyle{assumption}
\newtheorem{assumption}{Assumption}[section]

\DeclareMathOperator*{\cov}{cov}
\DeclareMathOperator*{\diag}{diag}
\DeclareMathOperator*{\trace}{trace}

\input cyracc.def

\begin{document}

\title[Moderate deviation principle]
{Moderate deviation principle for exponentially ergodic Markov
chain}
\author{B. Delyon}
\address{Universit\'e de Rennes 1, IRISA, Campus de Beaulieu, 35042 Rennes Cedex,
France.} \email{bernard.delyon@univ-rennes1.fr}
\author{A. Juditsky}
\address{University Joseph Fourier of Grenoble, France}
\email{juditsky@inrialpes.fr}

\author{R. Liptser}
\address{Department of Electrical Engineering-Systems,
Tel Aviv University, 69978 Tel Aviv Israel}
\email{liptser@eng.tau.ac.il}

\maketitle
\begin{abstract}
For $\frac{1}{2}<\alpha<1$, we propose the MDP analysis for family
$$
S^\alpha_n=\frac{1}{n^\alpha}\sum_{i=1}^nH(X_{i-1}), \ n\ge 1,
$$
where $(X_n)_{n\ge 0}$ be a homogeneous ergodic Markov chain,
$X_n\in \mathbb{R}^d$,  when the spectrum of operator $P_x$ is
continuous. The vector-valued function $H$ is not assumed to be
bounded but the Lipschitz continuity of $H$ is required. The main
helpful tools in our approach are Poisson equation and Stochastic
Exponential; the first enables to replace the original family by
$\frac{1}{n^\alpha}M_n$ with a martingale $M_n$ while the second
to avoid the direct Laplace transform analysis.

\end{abstract}

\section{Introduction and discussion}
\label{sec-1} Let $(X_n)_{n\ge 0}$ be a homogeneous ergodic Markov
chain, $X_n\in \mathbb{R}^d$   with the transition probability
kernel for $n$ steps: $P^{(n)}_x=P^{(n)}(x,dy)$ (for brevity
$P^{(1)}_x:=P_x$) and the unique invariant measure $\mu$.

Let $H$ be a measurable function
$\mathbb{R}^d\stackrel{H}{\to}\mathbb{R}^p$ with
$\int_{\mathbb{R}^d}|H(z)|\mu (dz)<\infty$ and
\begin{equation}\label{nul}
\int_{\mathbb{R}^d}H(z)\mu (dz)=0.
\end{equation}
Set
$$
S^\alpha_n=\frac{1}{n^\alpha}\sum_{i=1}^nH(X_{i-1}), \ n\ge 1; \
(0.5<\alpha<1).
$$

In this paper, we examine the moderate deviation principle (MDP)
for families $(S^\alpha_n)_{n\ge 1}$ when the spectrum of operator
$P_x$ is continuous.

It is well known that for bounded $H$ satisfying \eqref{nul} ((H)
- condition), the most MDP compatible Markov chains are
characterized by {\it eigenvalues gap condition {\rm (EG)}} (see
Wu, \cite{Wu}, \cite{Wu1}, Gong and Wu, \cite{GongWu5}, and
citations therein):
\begin{quote}
{\it the unit is an isolated, simple and the only eigenvalue with
modulus 1 of the transition probability kernel $P_x$}.
\end{quote}
In the framework of (H)-(EG)  conditions, the MDP is valid with
the rate of speed $n^{-(2\alpha-1)}$ and the rate function
$I(y),y\in\mathbb{R}^d$
\begin{equation}\label{I(x)2}
I(y)=
\begin{cases}
\frac{1}{2}\|y\|^2_{B^\oplus}, & B^\oplus By=y
\\
\infty, & \text{otherwise},
\end{cases}
\end{equation}
where  $B^\oplus$ is the pseudoinverse matrix (in Moore-Penrose
sense, see e.g.\cite {Albert}) for the matrix
\begin{multline}\label{B}
B=\int_{\mathbb{R}^d}H(x)H^*(x)\mu(dx)
\\
+\sum_{n\ge 1}\int_{\mathbb{R}^d}
\left[H(x)(P^{(n)}_xH)^*+(P^{(n)}_xH)H^*(x)\right]\mu(dx)
\end{multline}
(henceforth, $^*$, $|\cdot|$,  and $\|\cdot\|_Q$ are the
transposition symbol, $\mathbb{L}^1$ norm and $\mathbb{L}^2$ norm
with the kernel $Q$ ($\|x\|_Q=\sqrt{\langle x,Qx\rangle}$)
respectively).

Owing to the quadratic form rate function, the MDP is an
attractive tool for an asymptotic analysis in many areas, say,
with thesis
\begin{center}
``MDP instead of CLT''
\end{center}

In this paper, we intend to apply the MDP analysis to Markov chain
defined by the recurrent equation
$$
X_n=f(X_{n-1},\xi_n), \ n\ge 1
$$
generated by i.i.d. sequence $(\xi_n)_{n\ge 1}$ of random vectors,
where $f$ is some vector-valued measurable function. Obviously,
the function $f$ and the distribution of $\xi_1$ might be
specified in this way $P_x$ satisfies (EG). For instance, if $d=1$
and
$$
X_n=f(X_{n-1})+\xi_n,
$$
then for bounded $f$ and Laplacian random variable $\xi_1$ (EG)
holds. However, (EG) fails for many useful in applications ergodic
Markov chains. For $d=1$, a typical example is Gaussian Markov
chain defined by the linear recurrent equation governed by i.i.d.
sequence of $(0,1)$-Gaussian random variables(here $|a|<1$)
$$
X_n=aX_{n-1}+\xi_n.
$$

\medskip
In this paper, we avoid a verification of (EG). Although our
approach is close to conceptions of ``Multiplicative Ergodicity''
(see Balaji  and Myen \cite{BalMyen}) and ``Geometrical
Ergodicity'' (see Kontoyiannis and Meyn, \cite{KuMeyn} and Meyn
and Tweedie, \cite{metw}),  Chen and Guillin,
\cite{ChenGui}) we do not follow explicitly these
methodologies.

Our main tools are the Poisson equation and the Puhalskii theorem
from \cite{P3}. The Poisson equation permits to reduce the MDP
verification for $(S^\alpha_n)_{n\ge 1}$ to
$(\frac{1}{n^\alpha}M_n)_{n\ge 1}$, where $M_n$ is a martingale
generated by Markov chain, while the Puhalskii theorem allows to
replace an asymptotic analysis for the Laplace transform of
$\frac{1}{n^\alpha}M_n$ by the asymptotic analysis for, so called,
{\it Stochastic Exponential}
\begin{equation}\label{StEx}
\mathscr{E}_n(\lambda)=\prod_{i=1}^nE\Big(\exp\Big[\Big\langle\lambda,
\frac{1}{n^\alpha} (M_i-M_{i-1})\Big\rangle\Big]\Big|X_{i-1}\Big),
\ \lambda\in\mathbb{R}^d
\end{equation}
being the product of the conditional Laplace transforms for
martingale increments.

An effectiveness of  the Poisson equation approach (method of
corrector) combined with the stochastic exponential is well known
from the proofs of functional central limit theorem (FCLT)  for
the family $(S^{0.5}_n)_{n\ge 1}$ (see, e.g. Papanicolaou, Stroock
and Varadhan \cite{PSV}, Ethier and Kurtz \cite{EK}, Bhattacharya
\cite{Br}, Pardoux and Veretennikov \cite{PV}) and with
$0.5<\alpha<1$ for the MDP analysis for the continuous time case
(see e.g. \cite{LipSpo}, \cite{LSV}).

\section{Formulation of main result}
\label{sec-2} We consider Markov chain $(X_n)_{n\ge 0}$,
$X_n\in\mathbb{R}^d$ defined by a nonlinear recurrent equation
\begin{equation}\label{4.1a}
X_n=f(X_{n-1},\xi_n),
\end{equation}
where $f=f(z,v)$ is a vector function with entries
$f_1(z,v),\ldots,f_d(z,v)$, $u\in \mathbb{R}^d$, $v\in
\mathbb{R}^p$ and $(\xi_n)_{n\ge 1}$ is i.i.d. sequence of random
vectors of the size $p$.

We fix the following assumptions.

\medskip
\begin{assumption}\label{N.1}
Entries of $f$ are Lipschitz continuous functions in the following
sense:
\begin{multline*}
|f_i(z_1\ldots,z_{j-1},z'_j,z_{j+1}\ldots,z_d,v_1,\ldots,v_p)
\\
- f_i(z_1\ldots,z_{j-1},z''_j,z_{j+1}\ldots,z_d,v_1,\ldots,v_p)|
\\
\le\varrho_{ij}|z'_j-z''_j|,
\\ \\
|f(z,v')-f(z,v'')|\le \varrho|v'-v''|,\\
\end{multline*}
where
$$
\max_{i,j}\varrho_{ij}=\varrho<1.
$$
\end{assumption}

\begin{assumption}\label{L.2}
For sufficiently small positive $\delta$,
$$
Ee^{\delta|\xi_1|}<\infty.
$$
\end{assumption}

\medskip
\noindent

\begin{theorem}\label{theo-4.1}
Under Assumptions \ref{N.1} and \ref{L.2}, the Markov chain is
ergodic with the invariant measure $\mu$ such that
$\int_{\mathbb{R}^d}|z|\mu(dz)<\infty$.  For any Lipschitz
continuous function $H$ with $\int_{\mathbb{R}^d}H(z)\mu(dz)=0$,
the family $(S^\alpha_n)_{n\ge 1}$ obeys the MDP in the metric
space $(\mathbb{R}^d,r)$ ($r$ is the Euclidean metric) with the
rate of speed $n^{-(2\alpha-1)}$ and the rate function given in
\eqref{I(x)2}.
\end{theorem}
\begin{remark}
The assumptions of Theorem \ref{theo-4.1} do not guarantee {\rm (EG)}.
The function $H$ satisfies the linear growth condition. Moreover,
an existence of the continuous component of $\xi_1$-distribution,
as in is not to be assumed.
\end{remark}

\bigskip
Consider now a linear recurrent equation
$$
X_n=AX_{n-1}+\xi_n
$$
governed by i.i.d. sequence  $(\xi_n)_{n\ge 1}$ of random vectors,
where $A$ is matrix of the size $d\times d$ and entries $A_{ij}$.
Now, Assumption \ref{N.1}  reads as: $\max_{ij}|A_{ij}|<1$. This
assumption is very restrictive and we replace it by more natural
one

\begin{assumption}\label{L.1}
The eigenvalues of $A$ lie within the unit circle.
\end{assumption}

\begin{theorem}\label{theo-1}
Under Assumption \ref{L.1}, the Markov chain is ergodic with the
invariant measure $\mu$ such that
$\int_{\mathbb{R}^d}\|z\|^2\mu(dz)<\infty$. For any Lipschitz
continuous function $H$ with
$\int_{\mathbb{R}^d}H(z)\mu(dz)<\infty$, the family
$(S^\alpha_n)_{n\ge 1}$ obeys the MDP in the metric space
$(\mathbb{R}^d,r)$ with the rate of speed $n^{-(2\alpha-1)}$ and
the rate function given in \eqref{I(x)2}.
\end{theorem}

\section{Preliminaries}
\label{sec-3}
\subsection{(EG)-(H) conditions}
\mbox{} \label{sec-3.1} To clarify our approach to the MDP
analysis, let us demonstrate its application for (EG)-(H) setting.

The (EG) condition provides the geometric ergodicity   of
$P^{(n)}_x$ to the invariant measure $\mu$ in the total variation
norm: there exist constants $K>0$ and $\varrho\in(0,1)$ such that
for any $x\in\mathbb{R}^d$
$$
\|P^{(n)}_x-\mu\|_{\sf tv}\le K\varrho^n, \ n\ge 1.
$$
So, (EG)-(H) conditions provide the existence of bounded function
\begin{equation}\label{ux}
U(x)=H(x)+\sum_{n\ge 1}P^{(n)}_xH
\end{equation}
which solves the Poisson equation
\begin{equation}\label{Poisson1}
H(x)=H(x)+P_xU.
\end{equation}
In view of the Markov property, $ \zeta_i:=U(X_i)-P_{X_{i-1}}U, \
i\ge 1 $ is the sequence of bounded martingale-differences with
respect to the filtration generated by Markov chain. Hence,
$M_n=\sum_{i=1}^n\zeta_i$ is the martingale with bounded
increments. With the help of Poisson's equation we get the
following decomposition
\begin{equation}\label{PoissDec}
\sum_{i=1}^nH(X_{i-1})=\underbrace{U(x)-U(X_n)}_{\rm
corrector}+M_n.
\end{equation}
Owing to the boundedness of $U$, the families $S^\alpha_n$ and
$\frac{1}{n^\alpha}M_n$ share the same MDP. This fact enables us
to verify the MDP for $(\frac{1}{n^\alpha}M_n)_{n\ge 1}$.

\medskip
Assume for a moment that $\zeta_i$'s are also are independent and
identically distributed random vectors. Recall that $E\zeta_1=0$
and denote $B=E\zeta_1\zeta^*_1$. Introduce  the Laplace transform
for $\frac{1}{n^\alpha}M_n$:
\begin{equation}\label{Lap0}
\mathscr{E}_n(\lambda)=\Big(Ee^{\langle\lambda,
\frac{\zeta_1}{n^\alpha}\rangle}\Big)^n, \lambda\in\mathbb{R}^d.
\end{equation}
It is well known that the MDP for $\frac{1}{n^\alpha}M_n$ is
provided by the following conditions: $B$ is not singular matrix
and
$$
\lim_{n\to\infty}n^{2\alpha-1}
\log\mathscr{E}_n(\lambda)=\frac{1}{2}\langle
\lambda,B\lambda\rangle, \ \lambda\in \mathbb{R}^d.
$$

\bigskip
The framework of this proof might be adapted to the case
considered in the paper. Instead of $B$, we introduce the
conditional covariance matrix $ E(\zeta_i\zeta^*_i|X_{i-1}) $ and,
instead of the Laplace transform \eqref{Lap0}, the stochastic
exponential \eqref{StEx}, having a form
$$
\mathscr{E}_n(\lambda)=\prod_{i=1}^nE\Big(e^{\langle\lambda,
\frac{\zeta_i}{n^\alpha}\rangle}\big|X_{i-1}\Big), \ \lambda\in
\mathbb{R},
$$
which is not the Laplace transform itself. The homogeneity of
Markov chain and the definition of $\zeta_i$ provide $
E(\zeta_i\zeta^*_i|X_{i-1})=B(X_{i-1}) $ for the matrix-valued
function
\begin{equation}\label{B(x)}
B(x)=P_xUU^*-P_xU\big(P_xU)^*.
\end{equation}
The Poisson equation \eqref{Poisson1} and its solution \eqref{ux}
allow to transform \eqref{B(x)} into
$$
B(x)=H(x)H^*(x)+\sum_{n\ge 1}\Big[H(x)\big(P^{(n)}_xH\big)^*+
\big(P^{(n)}_xH\big) H^* \Big],
$$
that is for $B$ defined in \eqref{B}
$$
\int_{\mathbb{R}^d}B(z)\mu(dz)=B.
$$

\medskip
Now, we are in the position to formulate

\smallskip
\noindent {\bf Puhalskii Theorem.} {\it Assume $B$ from \eqref{B}
is nonsingular matrix and for any $\varepsilon>0$,
$\lambda\in\mathbb{R}^d$
\begin{equation}\label{2.7Puh11}
\lim_{n\to\infty}\frac{1}{n^{2\alpha-1}}\log
P\Big(\Big|n^{2\alpha-1}
\log\mathscr{E}_n(\lambda)-\frac{1}{2}\langle\lambda,B\lambda\rangle\Big|>
\varepsilon\Big)=-\infty.
\end{equation}

Then, the family $\frac{1}{n^\alpha}M_n$, $n\ge 1$ possesses the
MDP in the metric space $(\mathbb{R}^d,r)$ ($r$ is the Euclidean
metric) with the rate of speed $n^{-(2\alpha-1)}$ and rate
function $I(y)=\frac{1}{2}\|y\|^2_{B^{-1}}$. }

\begin{remark}
The condition \eqref{2.7Puh11} is provided by
\begin{equation}\label{I1}
\begin{aligned}
& \lim_{n\to\infty}\frac{1}{n^{2\alpha-1}}\log P
\Big(\frac{1}{n}\Big|\sum_{i=1}^n\Big\langle\lambda,\big[B(X_{i-1})-B\big]\lambda
\Big\rangle \Big|>\varepsilon\Big)=-\infty
\\
& \lim_{n\to\infty}\frac{1}{n^{2\alpha-1}}\log P
\Big(\frac{1}{6n^{1+\alpha}}\sum_{i=1}^n
E\Big[|\zeta_i|^3e^{n^{-\alpha}|\zeta_i|}\big|X_{i-1}\Big]
>\varepsilon\Big)=-\infty.
\end{aligned}
\end{equation}
\end{remark}

The second condition in \eqref{I1} is implied by the boundedness
of $|\zeta_i|$'s. The first part in \eqref{I1} is known as Dembo's
conditions, \cite{D}, formulated as follows: for any
$\varepsilon>0$, $\lambda\in\mathbb{R}^d$
$$
\varlimsup_{n\to\infty}\frac{1}{n}\log P
\Big(\frac{1}{n}\Big|\sum_{i=1}^n\Big\langle\lambda,\big[B(X_{i-1})-B\big]\lambda
\Big\rangle \Big|>\varepsilon\Big)<0.
$$
In order to verify the first condition in \eqref{I1}, we will
follow to the framework of Poisson's equation technique. We
introduce the function
$$
h(x)=\big\langle\lambda,\big[B(x)-B\big]\lambda\big\rangle
$$
which is bounded and
$$
\int_{\mathbb{R}^d}h(z)\mu(dz)=0.
$$
Then the function $ u(x)=h(x)+\sum_{n\ge 1}P^{(n)}_xh $ is well
defined and solves the Poisson equation $ u(x)=h(x)+P_xu. $
Similarly to \eqref{PoissDec}, we have
$$
\frac{1}{n}\sum_{i=1}^nh(X_{i-1})=\frac{u(x)-u(X_n)}{n}+
\frac{m_n}{n},
$$
where $m_n=\sum_{i=1}^nz_i$ is the martingale with bounded
martingale-differences $(z_i)_{i\ge 1}$. Since $u$ is bounded, the
first condition in \eqref{I1} is reduced to
\begin{equation}\label{3.13a}
\lim_{n\to\infty}\frac{1}{n^{2\alpha-1}}\log
P\big(|m_n|>n\varepsilon\big)=-\infty
\end{equation}
while \eqref{3.13a} is provided by Theorem \ref{theo-A.1} in
Appendix which states that \eqref{3.13a} holds for any martingale
with bounded increments.

\subsubsection{Singular $B$}
\label{sec-3.1.1} Though to the conditions from \eqref{I1} remain
to hold when $B$ is singular, the Puhalskii theorem is no longer
valid. Nevertheless, we shall use this theorem as an auxiliary
tool.

It is well known that the family $\frac{M_n}{n^\alpha}$, $n\ge 1$
obeys the MDP with the rate of speed $n^{-(2\alpha-1)}$ and the
rate function given in \eqref{I(x)2} provided that
\begin{equation}\label{2.21sing1}
\begin{aligned}
& \varlimsup_{C\to \infty}\varlimsup_{n\to
\infty}\frac{1}{n^{2\alpha-1}}P
\Big(\Big\|\frac{M_n}{n^\alpha}\Big\|>C\Big)=-\infty
\\
& \varlimsup_{\varepsilon\to 0}\varlimsup_{n\to
\infty}\frac{1}{n^{2\alpha-1}}P
\Big(\Big\|\frac{M_n}{n^\alpha}-y\Big\|\le\varepsilon\Big)\le
-I(y)
\\
& \varliminf_{\varepsilon\to 0}\varliminf_{n\to
\infty}\frac{1}{n^{2\alpha-1}}P
\Big(\Big\|\frac{M_n}{n^\alpha}-y\Big\|\le \varepsilon\Big)\ge
-I(y).
\end{aligned}
\end{equation}
The first condition in \eqref{2.21sing1} provides the exponential
tightness in the metric $r$ while the next others the local MDP.

In order to verify of \eqref{2.21sing1}, we introduce
``regularized'' family $\frac{M^\beta _n}{n^\alpha},n\ge 1$ with
$$
M^\beta _n=M_n+\sqrt{\beta }\sum_{i=1}^n\theta_i,
$$
where $\beta $ is a positive parameter and $(\theta_i)_{i\ge 1}$
is a sequence of zero mean i.i.d. Gaussian random vectors with
$\cov(\theta_1,\theta_1)=:\mathbf{I}$ ($\mathbf{I}$ is the unite
matrix). The Markov chain and  $(\theta_i)_{i\ge 1}$ are assumed
to be independent objects.

It is clear that for this setting the matrix $B$ is transformed
into the positive definite matrix $B_\beta=B+\beta\mathbf{I}$.
Now, the Puhalskii theorem is applicable and guarantees the MDP
with the same rate of speed and the rate function
$I_\beta(y)=\frac{1}{2}\|y\|^2_{B^{-1}_ \beta }$. As the corollary
of the MDP, $\frac{M^\beta _n}{n^\alpha},n\ge 1$ is exponentially
tight (see \cite{P1}) and  obeys the local MDP:
\begin{equation}\label{2.21song1}
\begin{aligned}
& \varlimsup_{C\to \infty}\varlimsup_{n\to
\infty}\frac{1}{n^{2\alpha-1}}P \Big(\Big\|\frac{M^\beta
_n}{n^\alpha}\Big\|>C\Big)=-\infty
\\
& \varlimsup_{\varepsilon\to 0}\varlimsup_{n\to
\infty}\frac{1}{n^{2\alpha-1}}P \Big(\Big\|\frac{M^\beta
_n}{n^\alpha}-u\Big\|\le\varepsilon\Big)\le -I_\beta (y)
\\
& \varliminf_{\varepsilon\to 0}\varliminf_{n\to
\infty}\frac{1}{n^{2\alpha-1}}P \Big(\Big\|\frac{M^\beta
_n}{n^\alpha}-u\Big\|\le\varepsilon\Big)\ge -I_\beta (y).
\end{aligned}
\end{equation}
Notice now that \eqref{2.21sing1} is implied by \eqref{2.21song1}
provided that
\begin{equation}\label{2.23+1}
\lim_{\beta \to 0}I_\beta (y)=
\begin{cases}
\frac{1}{2}\|u\|^2_{B^\oplus}, & B^\oplus By=y
\\
\infty, & \text{otherwise}
\end{cases}
\end{equation}
and
\begin{equation}\label{2.24end1}
\lim_{\beta \to 0}\varlimsup_{n\to \infty}\frac{1}{n^{2\alpha-1}}P
\Big(\Big\|\frac{\sqrt{\beta
}}{n^\alpha}\sum_{i=1}^n\theta_i\Big\|> \eta\Big)=-\infty, \quad \
\forall \ \eta>0.
\end{equation}

Let $T$ be the orthogonal matrix transforming $B$ to the diagonal
form: $ \diag(B)=T^*BT. $ Then, owing to
$$
2I_\beta (y)=y^*(\beta  I+B)^{-1}y=y^*T(\beta
I+\diag(B))^{-1}T^*y,
$$
for $y=B^\oplus By$ we have (recall that $B^\oplus
BB^\oplus=B^\oplus$, see \cite{Albert})
$$
\begin{aligned}
2I_\beta (y)&=y^*B^\oplus BT(\beta I+\diag(B))^{-1}T^*y
\\
&=y^*B^\oplus TT^* BT(\beta I+\diag(B))^{-1}T^*y
\\
&=y^*B^\oplus T\diag(B)(\beta I+\diag(B))^{-1}T^*y
\\
&\xrightarrow[\beta \to 0]{}y^*B^\oplus
T\diag(B)\diag^\oplus(B)T^*y
\\
&=y^*B^\oplus T\diag(B)T^*T\diag^\oplus(B)T^*y
\\
&=y^*B^\oplus BB^\oplus u=u^*B^\oplus y=\|y\|^2_{B\oplus}=2I(y).
\end{aligned}
$$

If $y\ne B^\oplus By$, $\lim_{\beta \to 0}2I_\beta (y)=\infty$.

Thus, \eqref{2.23+1} holds true.

\medskip
Since $(\theta_i)_{i\ge 1}$ is i.i.d. sequence of random vectors
and entries of $\xi_1$ are i.i.d. $(0,1)$-Gaussian random
variables, the verification of \eqref{2.24end1} is reduced to
\begin{equation}\label{Ibeta1}
\lim_{\beta \to
0}\varlimsup_{n\to\infty}\frac{1}{n^{2\alpha-1}}\log P\Big(
\Big|\sum_{i=1}^n\vartheta_i\Big|>\frac{n^\alpha\eta}{\sqrt{\beta
}}\Big) =-\infty,
\end{equation}
where $(\vartheta_i)_{i\ge 1}$ is a sequence of i.i.d.
$(0,1)$-Gaussian random variables. It is clear that \eqref{Ibeta1}
is equivalent to
$$
\lim_{\beta \to
0}\varlimsup_{n\to\infty}\frac{1}{n^{2\alpha-1}}\log P\Big(
\pm\sum_{i=1}^n\vartheta_i>\frac{n^\alpha\eta}{\sqrt{\beta }}\Big)
=-\infty
$$
and, moreover, it suffices to establish ``+'' only.

By the Chernoff inequality with $\lambda>0$, we find that
$$
P\big(\sum_{i=1}^n\theta_i>\frac{n^\alpha\eta}{\sqrt{\beta }}\big)
\le \exp\Big(-\lambda\frac{n^\alpha\eta}{\sqrt{\beta
}}+n\frac{\lambda^2} {2}\Big)
$$
while the choice of $\lambda=\frac{n^\alpha\eta}{n\sqrt{\beta }}$
provides
$$
\frac{1}{n^{2\alpha-1}}\log P\Big(
\sum_{i=1}^n\eta_i>\frac{n^\alpha\eta}{\sqrt{\beta }}\Big)\le -
\frac{\eta^2}{2\beta }\xrightarrow[\beta \to 0]{}-\infty.
$$

\subsection{Virtual scenario}
\label{sec-3.2} \mbox{}

- (EG)-(H) are not assumed

- the ergodicity of Markov is checked

- $H$ is chosen to hold \eqref{nul}.

\medskip
\noindent {\bf (1)} We assume \eqref{ux}. Then the function $U$
solves the Poisson equation and the decomposition from
\eqref{PoissDec} is valid with $M_n=\sum_{i=1}^n\zeta_i$, where
$\zeta_i=u(X_i)-P_{X_{i-1}}u$. We assume that
$$
\begin{aligned}
& E\zeta^*_i\zeta_i\le{\rm const.}
\\
&
E\Big[|\zeta_i|^3e^{n^{-\alpha}|\zeta_i|}\big|X_{i-1}\Big]\le{\rm
const.}
\end{aligned}
$$

\medskip
\noindent {\bf (2)} With $B(x)$ and $B$ are defined in
\eqref{B(x)} and \eqref{B} respectively, set
$$
h(x)=\Big\langle \lambda,\big[B(x)-B\big]\lambda\Big\rangle, \
\lambda\in\mathbb{R}^d.
$$

We assume that

(i) $u(x)=h(x)+\sum_{n\ge 1}P^{(n)}_xh$ is well defined

(ii) for $z_i=u(X_i)-P_{X_{i-1}}u$,
$$
\begin{aligned}
& Ez^2_i\le{\rm const.}
\\
& E\Big[|z_i|^3e^{n^{-\alpha}|z_i|}\big|X_{i-1}\Big]\le{\rm
const.}
\end{aligned}
$$

\medskip
\noindent {\bf (3)} We assume that for any $\varepsilon>0$
\begin{eqnarray*}
& \lim\limits_{n\to\infty}\frac{1}{n^{2\alpha-1}} \log
P\big(|U(X_n)|>n^\alpha\varepsilon\big)=-\infty
\\
& \lim\limits_{n\to\infty}\frac{1}{n^{2\alpha-1}} \log
P\big(|u(X_n)|>n^\alpha\varepsilon\big)=-\infty.
\end{eqnarray*}

\medskip
Notice that (EG)-(H) provide {\bf (1)}-{\bf (3)} and even if
(EG)-(H) fail, {\bf (1)}-{\bf (3)} may fulfill. Moreover, {\bf
(1)}-{\bf (3)} guarantee the validity for all steps of the proof
given in Section \ref{sec-3.1}, so that if the ergodic property of
Markov chain hold, {\bf (1)}-{\bf (3)} provides the MDP.

\medskip
We use this scenario for the proofs of Theorems \ref{theo-4.1} and
\ref{theo-1}.

\section{The proof of Theorem \ref{theo-4.1}}
\label{sec-4}

Here, we follow the virtual scenario from Section \ref{sec-3.2}.

\subsection{Ergodic property}

\begin{lemma}\label{lem-4.1}
Under Assumption \ref{N.1}, $(X_n)_{n\ge 0}$ possesses the unique
probability invariant measure $\mu$ with
$\int_{\mathbb{R}^d}|z|\mu(dz)<\infty$.
\end{lemma}

\begin{proof}
Let us initialize the recursion given in \eqref{4.1a} by a random
vector $X_0$, independent of $(\xi_n)_{n\ge 1}$, with the
distribution function $\nu$ such that
$\int_{\mathbb{R}^d}|x|\nu(dx)<\infty$ and denote $
\mu^n(dz)=\int_{\mathbb{R}^d}P^{(n)}_x(dz)\nu(dx). $ We show that
the family $(\mu^n)_{n\ge 1}$ is tight in the Levy-Prohorov
metric:
$$
\lim\limits_{k\to\infty}\varlimsup\limits_{n\to\infty}\mu^n(|z|>k)=0.
$$
Since by the Chebyshev inequality $\mu^n(|z|>k)\le
\frac{E|X_n|}{k}$, it suffices to show that
\begin{equation}\label{m4.2}
\sup_{n\ge 1}E|X_n|<\infty.
\end{equation}
By Assumption \ref{N.1}, we have
$$
\begin{aligned}
|X_n|&=|f(0,\xi_n)+(f(X_{n-1},\xi_n)-f(0,\xi_n))|
\\
&\le |f(0,\xi_n)|+|f(X_{n-1},\xi_n)-f(0,\xi_n))|
\\
&\le |f(0,\xi_n)|+\varrho|X_{n-1}|
\\
&\le |f(0,0)|+\ell|\xi_n|+\varrho|X_{n-1}|.
\end{aligned}
$$
Hence, $E|X_n|\le |f(0,0)|+\ell E|\xi_1|+\varrho E|X_{n-1}|$.
Since $E|X_0|<\infty$, we have $E|X_n|<\infty$ and, moreover,
\begin{equation}\label{m4.3}
E|X_n|\le E|X_0|+\frac{|f(0,0)|+\ell E|\xi_1|} {1-\varrho}, \quad
n\ge 1,
\end{equation}
i.e. \eqref{m4.2} holds true.

\smallskip
Further, by the Prohorov theorem, the sequence of $\mu^n$,
$n\nearrow\infty$, contains further subsequence $\mu^{n'}$,
$n'\nearrow\infty$ converging in the Levy-Prohorov metric to a
limit $\mu$ being the probability measure on $\mathbb{R}^d$. In
other words, for any bounded and continuous function $g$
$$
\lim_{n'\to\infty}\int_{\mathbb{R}^d}g(z)\mu^{n'}(dz)=\int_{\mathbb{R}^d}g(z)\mu(dz).
$$
With $L>0$, taking $g(z)=L\wedge |z|$, we obtain
$$
\int_{\mathbb{R}^d}(L\wedge|z|)\mu(dz)=\lim_{n'\to\infty}E(L\wedge|X_{n'}|)
\le\varlimsup_{n\to\infty}E|X_n|\le\text{const.}
$$
(see \eqref{m4.3}). Then, by the monotone convergence theorem, it
holds that
$$
\int_{\mathbb{R}^d}|z|\mu(dz)\le\varlimsup_{n\to\infty}E|X_n|<\infty.
$$

Now, we show that $\mu$ is an invariant measure of the Markov
chain, that is for any nonnegative, bounded and measurable
function $g$
\begin{equation}\label{4.1aa}
\int_{\mathbb{R}^d}g(x)\mu(dx)=\int_{\mathbb{R}^d}P_xg\mu(dx).
\end{equation}
It suffices to verify \eqref{4.1aa} for continuous function only.
The general statement is obtained by a monotonic type
approximation.

For notational convenience, write $X^x_n$ and $X^\nu_n$, if
$X_0=x$ and $X_0$ is distributed in the accordance with $\nu$.
Making use Assumption \ref{N.1}, we find that
$$
|X^x_n-X^\nu_n|\le\varrho|X^x_{n-1}-X^\nu_{n-1}|, \ n\ge 1
$$
and claim that $|X^x_n-X^\nu_n|$ converges to zero exponentially
fast as long as $n\to\infty$. For any $x\in\mathbb{R}^d$, the
latter provides
$$
\lim_{n'\to\infty}Eg(X^x_{n'})=\int_{\mathbb{R}^d}g(x)\mu(dx).
$$
This and the fact that $(X_n)_{n\ge 0}$ is the homogeneous Markov
chain also imply
$$
\lim_{n'\to\infty}Eg(X^x_{n'+1})=\int_{\mathbb{R}^d}g(z)\mu(dz).
$$
On the other hand, taking into the consideration
$Eg(X^x_{n'+1})=EP_{X^x_{n'}}g$, the above relation is nothing but
$
\lim\limits_{n'\to\infty}EP_{X^x_{n'}}g=\int_{\mathbb{R}^d}g(z)\mu(dz).
$ The next key tool in the proof is the Feller property: for any
bounded and continuous $g$, the function $ P_xg=Eg(f(x,\xi_1)) $
is the continuous as well. So, by the Feller property $
\lim\limits_{n'\to\infty}EP_{X^x_{n'}}g=\int_{\mathbb{R}^d}P_xg\mu(dx).
$

Thus, \eqref{4.1aa} holds.

Assume $\mu'$ is another invariant probability measure, $\mu'\ne
\mu$. Then taking two random vectors $X^\mu_0$ and $X^{\mu'}_0$,
distributed in the accordance to $\mu$ and $\mu'$ respectively and
independent of $(\xi_n)_{n\ge 1}$, we create two stationary Markov
chains $(X^\mu_n)$ and $(X^{\mu'}_n)$ defined on the same
probability space as:
$$
\begin{aligned}
X^\mu_n=f(X^\mu_{n-1},\xi_n)
\\
X^{\mu'}_n=f(X^{\mu'}_{n-1},\xi_n).
\end{aligned}
$$
By Assumption \ref{N.1},
$|X^\mu_n-X^{\mu'}_n|\le\varrho|X^\mu_{n-1}-X^{\mu'}_{n-1}|$, i.e.
$\lim\limits_{n\to\infty}|X^\mu_n-X^{\mu'}_n|=0$ what contradicts
$\mu\neq \mu'$.
\end{proof}
\subsection{The verification of (1)}
\label{sec-4.2} Let $K$ be the Lipschitz constant for $H$. Then
$|H(x)|\le |H(0)|+K|x|$ and
$\int_{\mathbb{R}^d}|H(z)|\mu(dz)<\infty.$ Hence, \eqref{nul} is
the correct assumption and for the stationary Markov chain
$X^\mu_n$ we have $EH(X^\mu_n)\equiv 0$. Then,
$$
\begin{aligned}
|EH(X^x_n)|&=|E(H(X^x_n)-H(X^\mu_n)|
\\
&\le K\varrho^nE|x-X^\mu_n|\le K(1+|x|)\varrho^n.
\end{aligned}
$$
Therefore, $\sum_{n\ge 1}|EH(X^x_n|\le\frac{K}{1-\varrho}(1+|x|)$.
Consequently, the function $U(x)$, given in \eqref{ux}, is well
defined and solves the Poisson equation.

Recall that $\zeta_i=U(X_i)-P_{X_{i-1}}U$.

\begin{lemma}\label{lem-2.1a}
The function $U(x)$ possesses the following properties:

{\rm 1)} $U(x)$ is Lipschitz continuous;

{\rm 2)} $P_x(UU^*)-P_xU(P_xU)^*$ is bounded and Lipschitz
continuous;

{\rm 3)} For sufficiently small $\delta>0$ and any $i\ge 1$
$$
E\Big(\big|U(X_i)-P_{\mbox{}_{X_{i-1}}}U\big|^3e^{\delta|U(X_i)-
P_{\mbox{}_{X_{i-1}}}U|}\big|X_{i-1}\Big)\le \text{\rm const.}
$$
\end{lemma}
\begin{proof}
1) Since  by Assumption \ref{N.1},
$$
|X^{x'}_n-X^{x''}_n|\le \varrho|X^{x'}_{n-1}-X^{x''}_{n-1}|, \
|X^{x'}_0-X^{x''}_0|\le |x'-x''|,
$$
we have
\begin{equation}\label{4.6pi}
\begin{aligned}
|U(x')-U(x'')|&\le |H(x')-H(x'')|+\sum_{n\ge
1}E|H(X^{x'}_n)-H(X^{x''}_n)|
\\
&\le \frac{K}{1-\varrho}|x'-x''|.
\end{aligned}
\end{equation}

\medskip
2) Recall (see \eqref{B(x)})
$$
P_x(UU^*)-P_xU(P_xU)^*=B(x)
$$
and denote $B_{pq}(x)$, $p,q=1,\ldots,d$ the entries of matrix
$G(x)$. Also, denote by $U_p(x)$, $p=1,\ldots,d$ the entries of
$U(x)$. Since $B(x)$ is nonnegative definite matrix, it suffices
to show only that $B_{pp}(x)$'s are bounded functions. Denote
$F(z)$ the distribution function of $\xi_1$. Taking into the
consideration \eqref{4.6pi} and Assumption \ref{N.1}, we get
$$
\begin{aligned}
B_{pp}(x)&=E\Big(U_p\big(f(x,\xi_1)\big)-\int_{\mathbb{R}^d}U_p\big(f(x,z)\big)
dF(z)\Big)^2
\\
&\le
\frac{(K\ell)^2}{(1-\varrho)^2}E\Big|\int_{\mathbb{R}^d}|\xi_1-z|
dF(z)\Big|^2 \le
4\frac{(K\ell)^2}{(1-\varrho)^2}E|\xi_1|^2<\infty.
\end{aligned}
$$

The Lipschitz continuity of $B_{pq}(x)$ is proved similarly. Write
$$
B_{pq}(x')-B_{pq}(x'')=:ab-cd,
$$
where
$$
\begin{aligned}
&
a=E\Big(U_p\big(f(x',\xi_1)\big)-\int_{\mathbb{R}^d}U_q\big(f(x',z)\big)dF(z)\Big)
\\
&
b=E\Big(U_q\big(f(x',\xi_1)\big)-\int_{\mathbb{R}^d}U_q\big(f(x',z)dF(z)\Big)
\\
&
c=E\Big(U_p\big(f(x'',\xi_1)\big)-\int_{\mathbb{R}^d}U_q\big(f(x'',z)\big)dF(z)\Big)
\\
&
d=E\Big(U_q\big(f(x'',\xi_1)\big)-\int_{\mathbb{R}^d}U_q\big(f(x'',z)\big)dF(z)\Big).
\end{aligned}
$$
Now, applying $ab-cd=a(b-d)+d(a-c)$ and taking into account
\eqref{4.6pi} and Assumption \ref{N.1}, we find that $|a|,|d|\le
\frac{2K\ell}{1-\varrho}E|\xi_1| $ and so
$$
|B_{pq}(x')-B_{pq}(x'')|\le\frac{4K^2\ell\varrho}{(1-\varrho)^2}E|\xi_1||x'-x''|.
$$

\medskip
3) By \eqref{4.6pi} and Assumption \ref{N.1}
$$
|U(X_i)-P_{X_{i-1}}U|\le
\frac{K\ell}{1-\varrho}\big(E|\xi_1|+|\xi_i|\big).
$$
\end{proof}

\subsection{The verification of (2)}
\label{sec-4.3} The properties of $G(x)$ to be bounded and
Lipschitz continuous provide the same properties for
$$
h(x)=\big\langle \lambda,\big[B(x)-B\big]\lambda\big\rangle.
$$

Hence {\bf (2)} is provided by {\bf (1)}.

\subsection{The verification of (3)}
\label{sec-4.4} Since $U$ and $u$ are Lipschitz continuous, they
possess the linear growth condition, e.g. $ |U(x)|\le C(1+|x|), \
\exists C>0. $ So, {\bf (3)} is reduced to the verification of
\begin{equation}\label{2.10a}
\lim_{n\to\infty}\frac{1}{n^{2\alpha-1}}\log P\Big(\big|X_n\big|>
\varepsilon n^\alpha\Big)=-\infty, \ \varepsilon>0.
\end{equation}

Due to Assumption \ref{N.1}, we have
$$
\begin{aligned}
|X_n|&\le |f(X_{n-1},\xi_n)|\le |f(0,\xi_n)|+\varrho|X_{n-1}|
\\
&\le |f(0,0)| +\varrho|X_{n-1}|+\ell|\xi_n|.
\end{aligned}
$$
Iterating this inequality with $X_0=x$ we obtain
$$
\begin{aligned}
|X_n|&\le \varrho^n |x|+|f(0,0)|\sum_{j=1}^n\varrho^{n-j}+\ell
\sum_{j=1}^n\varrho^{n-j} |\xi_j|
\\
&\le
|x|+\frac{|f(0,0)|}{1-\varrho}+\ell\sum_{j=0}^{n-1}\varrho^j|\xi_{n-j}|.
\end{aligned}
$$
Hence, \eqref{2.10a} is reduced to
\begin{equation}\label{4.9m}
\lim_{n\to\infty}\frac{1}{n^{2\alpha-1}}\log
P\Big(\sum_{j=0}^{n-1}\varrho^j|\xi_{n-j}| \ge
n^\alpha\varepsilon\Big)=-\infty.
\end{equation}
We verify \eqref{4.9m} with the help of Chernoff's inequality:
with $\delta$, involving in Assumption \ref{L.2}, and
$\gamma=\frac{\delta}{1-\varrho}$
$$
\begin{aligned}
P\Big(\sum_{j=0}^{n-1}\varrho^j|\xi_{n-j}| \ge
n^\alpha\varepsilon\Big)\le e^{-n^\alpha
\gamma\varepsilon}Ee^{\sum_{j=0}^{n-1}\gamma\varrho^j|\xi_{n-j}|}.
\end{aligned}
$$
The i.i.d. property for $\xi_j$'s provides
$$
Ee^{\sum_{j=0}^{n-1}\gamma\varrho^j|\xi_{n-j}|}=
Ee^{\sum_{j=0}^{n-1}\gamma\varrho^j|\xi_1|}\le
Ee^{\sum_{j=0}^\infty\gamma\varrho^j|\xi_1|}=Ee^{\delta|\xi_1|}<\infty
$$
and we get
$$
\frac{1}{n^{2\alpha-1}}\log
P\Big(\sum_{j=0}^{n-1}\varrho^j|\xi_{n-j}| \ge
n^\alpha\varepsilon\Big) \le -n^{1-\alpha}\delta\varepsilon
+\frac{\log
Ee^{\delta|\xi_1|}}{n^{2\alpha-1}}\xrightarrow[n\to\infty]{}-\infty.
$$

\section{The proof of Theorem \ref{theo-1}}
\label{sec-5}

The proof of this theorem differs from the proof of Theorem
\ref{theo-4.1} only in some details concerning to (L.1). So, only
these parts of the proof are given below.

\subsection{Ergodic property and invariant measure}
\label{sec-5.1}

Introduce $(\widetilde{\xi}_n)_{n\ge 1}$ the independent copy of
$(\xi_n)_{n\ge 1}$. Owing to
$$
X_n=A^nx+\sum_{i=1}^nA^{n-i}\xi_i=A^nx+\sum_{i=0}^{n-1}A^i\xi_{n-i},
$$
we introduce
\begin{equation}\label{4.2x}
\widetilde{X}_n=A^nx+\sum_{i=0}^{n-1}A^i\widetilde{\xi}_i
\end{equation}
and notice that the i.i.d. property of $(\xi_i)_{i\ge 1}$ provides
$ (X_n)_{n\ge 0}\stackrel{\rm law}{=}(\widetilde{X}_n)_{n\ge 0}. $

By Assumption \ref{L.1}, $A^n\to 0$, $n\to\infty$, exponentially
fast. Particularly,
$$
\sum_{i=0}^\infty \trace\big(A^i\cov(\xi_1,\xi_1)
(A^i)^*\big)<\infty,
$$
so that $
\lim\limits_{n\to\infty}\widetilde{X}_n=\sum_{i=0}^\infty A^i
\widetilde{\xi}_i$ a.s. and in $\mathbb{L}^2$ norm.

Thus, the invariant measure $\mu$ is generated by the distribution
function of $\widetilde{X}_\infty$. In addition,
$E\|\widetilde{X}_\infty\|^2=\sum\limits_{i=0}^\infty
\trace\big(A^i\cov(\xi_1,\xi_1) (A^i)^*\big)$, so that
$$
\int_{\mathbb{R}^d}\|z\|^2\mu(dz)<\infty.
$$

\subsection{The verification of (1) and {\bf (2)}}
\label{sec-5.2}

Due to
$$
(X^{x'}_n-X^{x''}_n)=A(X^{x'}_{n-1}-X^{x''}_{n-1}),
$$
we have $(X^{x'}_n-X^{x''}_n)=A^n(x'-x'')$. Let us transform the
matrix $A$ into a Jordan form $ A = T J T^{-1} $ and notice that
$A^n=TJ^nT^{-1}$. It is well known that the maximal absolute value
of entries of $J^n$ is $n|\lambda|^n$, where $|\lambda|$ is the
maximal absolute value among eigenvalues of $A$. By Assumption
\ref{L.1}, $|\lambda|<1$. So, there exist $K>0$ and $\varrho<1$
such that $|\lambda|<\varrho$. Then, entries $A^n_{pq}$ of $A^n$
are evaluated as: $|A^n_{pq}|\le K\varrho^n$. Hence,
$|X^{x'}_n-X^{x''}_n|\le K\varrho^n|x'-x''|$, $n\ge 1$, and the
verification of {\bf (1), (2)} is in the framework of Section
\ref{sec-3}.

\subsection{The verification of (3)}
\label{sec-5.3} As in Section \ref{sec-3}, the verification of
this property is reduced to
\begin{equation}\label{2.10a}
\lim_{n\to\infty}\frac{1}{n^{2\alpha-1}}\log P\Big(\big|X_n\big|>
\varepsilon n^\alpha\Big)=-\infty, \ \varepsilon>0.
\end{equation}
In \eqref{2.10a}, we may replace $X_n$ by its copy
$\widetilde{X}_n$ defined in \eqref{4.2x}. Notice also that
$$
|\widetilde{X}_n|\le |A^nx|+\sum_{i=0}^\infty
\max_{pq}|A^i_{pq}||\widetilde{\xi}|.
$$
As was mentioned above, $|A^i_{pq}|\le K\varrho^j$ for some $K>0$
and $\varrho\in(0,1)$. Hence, it suffices to verify
$$
\lim_{n\to\infty}\frac{1}{n^{2\alpha-1}}\log
P\Big(\sum_{i=0}^\infty \varrho^i |\xi_i|>\varepsilon
n^\alpha\Big)=-\infty, \ \varepsilon>0
$$
what be going on similarly to corresponding part of the proof in
Section \ref{sec-3}.

\section{Exotic example}
\label{sec-2.3} Let $(X_n)_{n\ge 0}$, $X_n\in\mathbb{R}$ and
$X_0=x$, be Markov chain defined by the recurrent equation
\begin{equation}\label{exmp1}
X_n=X_{n-1}-m\frac{X_{n-1}}{|X_{n-1}|}+\xi_n
\end{equation}
where $m$ is a positive parameter, $(\xi_n)$ is i.i.d. sequence of
zero mean random variables with
$$
Ee^{\delta|\xi_1|}<\infty,  \ \text{for some $\delta>0$},
$$
and let $\frac{0}{0}=0$.

Although the virtual scenario is not completely verifiable here we
show that for
$$
H(x)=\frac{x}{|x|}
$$
the family $(S^\alpha_n)_{n\ge 1}$ possesses the MDP provided that
\begin{equation}\label{mdel}
m>\frac{1}{\delta}\log Ee^{\delta|\xi_1|}.
\end{equation}

Indeed, by \eqref{exmp1} we have
$$
\frac{1}{n^\alpha}\sum_{k=1}^n\frac{X_{k-1}}{|X_{k-1}|}=
\frac{1}{m}\frac{(X_n-x)}{n^\alpha}+\frac{1}{n^\alpha}
\sum_{k=1}^n\frac{\xi_k}{m}.
$$
The family
$\big(\frac{1}{n^\alpha}\sum\limits_{k=1}^n\frac{\xi_k}{m}\big)_{n\ge
1}$ possesses the MDP with the rate of speed $n^{-(2\alpha-1)}$
and the rate function $I(y)=\frac{m^2}{2E\xi^2_1}y^2$. Then, the
family $(S^\alpha_n)_{n\ge 1}$ obeys the same MDP provided that
$\big(\frac{X_n-x}{n^\alpha}\big)_{n\ge 1}$ is exponentially
negligible family with the rate $n^{-(2\alpha-1)}$. This
verification is reduced to
\begin{equation}\label{negl}
\lim_{n\to\infty}\frac{1}{n^{2\alpha-1}}\log
P\big(|X_n|>n^\alpha\varepsilon\big)= -\infty, \ \varepsilon>0.
\end{equation}
By the Chernoff inequality $
P\big(|X_n|>n^\alpha\varepsilon\big)\le e^{-\delta
n^\alpha\varepsilon} Ee^{\delta|X_n|}, $ that is \eqref{negl}
holds if $\sup\limits_{n\ge 1}Ee^{\delta|X_n|}<\infty$ for some
$\delta>0$. We show that the latter holds true for $\delta$
involved in \eqref{mdel}. A helpful tool for this verification is
the inequality $\big|z-m\frac{z}{|z|}\big|\le \big||z|-m\big|$.
Write
$$
\begin{aligned}
Ee^{\delta|X_n|}&=Ee^{\delta|X_n|}I(|X_{n-1}|\le
m)+Ee^{\delta|X_n|}I(|X_{n-1}|>m)
\\
&\le e^{\delta m}Ee^{\delta|\xi_1|}+e^{-\delta
m}Ee^{\delta|\xi_1|}Ee^{\delta|X_{n-1}|}.
\end{aligned}
$$
Set $\ell=e^{\delta m}Ee^{\delta|\xi_1|}$ and $\varrho=e^{-\delta
m}Ee^{\delta|\xi_1|}$. By \eqref{mdel}, $\varrho<1$. Hence,
$V(x)=e^{\delta|x|}$ is the Lyapunov function: $ P_xV\le\varrho
V(x)+\ell. $ Consequently,
$$
EV(X_n)\le \varrho EV(X_n)+\ell, \  n\ge 1
$$
and so, $\sup_{n\ge 1}EV(X_n)\le V(x)+\frac{\ell}{1-\varrho}$.

\appendix
\section{Exponentially integrable martingale-differences}

Let $\zeta_n=(\zeta_n)_{n\ge 1}$ be a martingale-difference  with
respect to some filtration $\mathscr{F}=(\mathscr{F}_n)_{n\ge 0}$
and $M_n=\sum\limits_{i=1}^n\zeta_i$ be the corresponding
martingale.

\begin{theorem}\label{theo-A.1}
Assume that for sufficiently small positive $\delta$ and any $i\ge
1$
\begin{equation}\label{A.1}
E\big(e^{\delta|\zeta_i|}|\mathscr{F}_{i-1}\big)\le {\rm const.}
\end{equation}
Then for any $\alpha\in (0.5,1)$
$$
\lim_{n\to\infty}\log \frac{1}{n^{2\alpha-1}}
P\big(|M_n|>n\varepsilon\big)=-\infty.
$$
\end{theorem}
\begin{proof}
It suffices to prove $
\lim\limits_{n\to\infty}\frac{1}{n^{2\alpha-1}}\log P\big(\pm
M'_n>n \varepsilon\Big)=-\infty $ and, moreover, it suffices to
verify ``+'' only  (``-'' is verified similarly).

\medskip
For fixed positive $\lambda$ and sufficiently large $n$, let us
introduce the stochastic exponential
$$
\mathscr{E}_n(\lambda)=\prod_{i=1}^nE\big(e^{\lambda
\frac{\zeta_i}{n}}\big|\mathscr{F}_{i-1}\big).
$$
A direct verification shows that
$$
E\exp\Big(\frac{\lambda M_n}{n}-\log
\mathscr{E}_n(\lambda)\Big)=1.
$$
We apply this equality for further ones
\begin{equation}\label{StTv}
\begin{aligned}
1&\ge EI\Big(M_n>n\varepsilon\Big) \exp\Big(\frac{\lambda
M_n}{n}-\log  \mathscr{E}_n(\lambda)\Big)
\\
&\ge EI\Big(M_n>n\varepsilon\Big) \exp\Big(\lambda\varepsilon-\log
\mathscr{E}_n(\lambda)\Big).
\end{aligned}
\end{equation}
Due to $E\big(\lambda\frac{\zeta_i}{n}|\mathscr{F}_{i-1}\big)=0$
and \eqref{A.1}, we find that
$$
\begin{aligned}
\log  \mathscr{E}_n(\lambda)&=\sum_{i=1}^n\log\Big(1+E\big[
e^{\lambda\frac{\zeta_i}{n}}-1-\lambda\frac{\zeta_i}{n}|\mathscr{F}_{i-1}\big]\Big)
\\
&
\le\sum_{i=1}^n\Big\{\frac{\lambda^2}{2n^2}E\big((\zeta_i)^2|X_{i-1}\big)+
\frac{\lambda^3}{6n^3}E\big(|\zeta_i|^3e^{\lambda\frac{|\zeta_i|}{n}}|\mathscr{F}
_{i-1}\big) \Big\}
\\
& \le K\Big[\frac{\lambda^2}{2n}+\frac{\lambda^3}{6n^2}\Big],
\end{aligned}
$$
where $K$ is some constant. This inequality, being incorporated
into \eqref{StTv}, provides
$$
1\ge EI\big(M_n>n\varepsilon\big)
\exp\Big(\lambda\varepsilon-K\Big[\frac{\lambda^2}{2n}+\frac{\lambda^3}{6n^2}
\Big]\Big).
$$
If $\varepsilon_0<3$, taking $\lambda=\varepsilon_0 nK^{-1}$, we
find that
$$
\frac{1}{n^{2\alpha-1}}\log P\big(M_n>n\varepsilon\big)\le
-\frac{\varepsilon^2n^{2(1-\alpha)}}{K}\Big(\frac{1}{2}-\frac{\varepsilon}{6}\Big)
\xrightarrow[n\to\infty]{}-\infty.
$$
Hence, for any verification $\varepsilon>0$ the desired statement
holds true.
\end{proof}

\end{document}